\documentclass[12pt]{article}
\usepackage{amsmath}
\usepackage{amssymb}
\usepackage{amsthm}
\usepackage{dsfont}
\usepackage{cite}

\def\nn{\nonumber}
\def\a{\alpha}

\def\g{\gamma}

\def\de{\delta}

\def\vk{\varkappa}

\def\la{\lambda}
\def\om{\omega}
\def\Om{\Omega}
\def\t{\theta}
\def\Up{\Upsilon}
\def\vp{\varphi}
\def\vt{\vartheta}
\def\ve{\varepsilon}
\def\ze{\zeta}
\def\wh{\widehat}
\def\wt{\widetilde}
\def\ov{\overline}

\def\BC{{\mathbb C}}
\def\BD{{\mathbb D}}
\def\BR{{\mathbb R}}

\def\cla{{\mathcal A}}
\def\clb{{\mathcal B}}
\def\clc{{\mathcal C}}

\def\clh{{\mathcal H}}

\def\clm{{\mathcal M}}
\def\cln{{\mathcal N}}

\def\clv{{\mathcal V}}

\def\mfa{{\mathfrak A}}
\def\ker{{\rm Ker}}

\newcommand{\E}{\mathrm{e}}
\newcommand{\I}{\mathrm{i}}
\def\mf{\mathfrak}

\newtheorem{Pa}{Paper}[section]
\newtheorem{Tm}[Pa]{{\bf Theorem}}
\newtheorem{La}[Pa]{{\bf Lemma}}

\newtheorem{Rk}[Pa]{{\bf Remark}}

\newtheorem{Dn}[Pa]{{\bf Definition}}
\newtheorem{Nn}[Pa]{{\bf Notation}}
\newtheorem{Pn}[Pa]{{\bf Proposition}}

\title{$S$-nodes, factorisation of spectral matrix functions and corresponding inequalities}

\author{Alexander Sakhnovich\footnote{This research    was supported by the
Austrian Science Fund (FWF) grant,  DOI: 10.55776/Y963.}}

\date{}
\begin{document}
\maketitle

\begin{flushleft}
Faculty of Mathematics, University of Vienna, \\
Oskar-Morgenstern-Platz 1, A-1090 Vienna, Austria\\
E-mail: oleksandr.sakhnovych@univie.ac.at
\end{flushleft}

\vspace{0.3em}

{\it To the memory of V.E. Katsnelson}

\vspace{0.3em} 

\begin{abstract} Using factorisation and Arov-Krein inequality results, we derive
important inequalities (in terms of $S$-nodes) in interpolation problems. 
\end{abstract}

\vspace{0.3em}


\vspace{0.3em}


\section{Introduction}\label{Intro}
\setcounter{equation}{0}
Let us consider some finite or infinite-dimensional Hilbert space $\clh$ and operators $A,S\in \clb(\clh)$, $\Pi \in \clb(\BC^{2p},\clh)$, which satisfy
the operator identity
\begin{align}& \label{I1}
AS-SA^*=\I \Pi J \Pi^*, \quad  J:=\begin{bmatrix} 0 & I_p \\ I_p & 0 \end{bmatrix},
\end{align}
where $\clb(\clh_1,\clh_2)$ stands for the set of bounded linear operators acting from the Hilbert space $\clh_1$ into Hilbert space $\clh_2$,
$\clb(\clh):=\clb(\clh,\clh)$, $\BC$ is (as usually) the complex plane, the symbol $\BC^{n\times p}$ denotes the set of $n\times p$ matrices 
with complex-valued entries ($\BC^n=\BC^{n\times 1}$), $\I$ stands for the  imaginary unit ($\I^2=-1$) and $I_p$ is the $p\times p$ identity matrix.
The operators $A^*$ and $\Pi^*$ in \eqref{I1} are adjoint to $A$ and $\Pi$, respectively, and $S$ is a self-adjoint operator ($S=S^*$).

{\it The triple $\{A,S,\Pi\}$ with the mentioned above properties is called  a self-adjoint $S$-node $($see \cite{SaSaR, SaL1, SaL2} and the references therein$)$}
or simply $S$-node because the $S$-nodes are always self-adjoint in this paper.

Using $S$-nodes, a wide class of interpolation problems is solved in a general way \cite{SaL2}.
For this purpose, $2p\times 2p$ matrix valued functions (so called {\it frames}) ${\mf A}(S,z)$ with the $p\times p$ blocks $\mfa_{ij}(S,z)=\mfa_{ij}(z)$
are introduced in \cite{SaL2} (see also the references therein):
\begin{align}& \label{I2}
\mfa(S,z)=\{\mfa_{ij}(z)\}_{i,j=1}^2=w_A(1/\ov{z})^*=I_{2p}-\I z\Pi^*(I-zA^*)^{-1}S^{-1}\Pi J,
\end{align}
where $I$ is the identity operator and
\begin{align}& \label{I3}
w_A(\la)= I_{2p}-\I J\Pi^*S^{-1}\big(A-\la I\big)^{-1}\Pi
\end{align}
is the transfer matrix function (matrix valued function) in Lev Sakhnovich form \cite{SaL1, SaL2}.  We partition the operator $\Pi$ into two blocks
\begin{align}& \label{I4}
\Pi = \begin{bmatrix} \Phi_1 & \Phi_2\end{bmatrix} \quad \big(\Phi_k\in \clb(\BC^{p},\clh)\big),
\end{align}
which generates the partition of $\mfa(S,z)$ into the blocks $\mfa_{ij}(z)$.

It is  assumed in important interpolation problems
(and will be usually assumed here) that 
\begin{align}& \label{I5}
S\geq \ve I \quad (\ve>0), \quad \ker \, \Phi_2=0,
\end{align} 
where $S\geq \ve I $ may be expressed in terms of the scalar products $(\cdot , \cdot )$ as $(Sf,f)\geq \ve (f,f)$ $(f \in\clh)$ and $\ker$ is
the kernel (nullspace) of the corresponding operator. (The strict operator inequality $S_1>S_2$ means that $\big((S_1-S_2)f,f\big)>0 $ for $f\not=0$.)
In view of \eqref{I5}, $S^{-1}$ exists and $S^{-1}\in \clb(\clh)$. Therefore, $\mfa(S,z)$ is well defined
in the points of invertibility of  $I-zA^*$. {\it We will assume that  the operators
$I-zA^*$ have  bounded inverse operators $($bounded inverses$)$ for $z\in \BC_+$, excluding, possibly, some isolated points, which are the poles of the
matrix function $\mfa(S,z)$.} Here, $\BC_+$ denotes open upper half-plane $\Im(z)>0$, where $\Im(z)$ is the imaginary part of $z$.

It easily follows from \eqref{I5} and identity \eqref{I1} (see also \cite{SaL1} or \cite[(1.88)]{SaSaR}) that
\begin{align}& \label{I6}
\mfa(S,z)J \mfa(S,\ov{\la})^*=J-\I (z-\la)\Pi^*(I-zA^*)^{-1}S^{-1}(I-\la A)^{-1}\Pi.
\end{align} 
It also follows from \eqref{I5} that $S^{-1}>0$. Hence, \eqref{I6} yields that $\break \mfa(S,z)J\mfa(S,z)^*\geq J$ for $\Im(z)>0$ ,
that is, for $z\in \BC_+$. 
Equivalently (see, e.g.,
 \cite[Corollary~E3]{SaSaR}) we have the inequality
\begin{align}& \label{I7}
\mfa(S,z)^*J\mfa(S,z)\geq J \qquad {\mathrm{for}} \qquad z\in \BC_+.
\end{align} 
Moreover, for the lower right block of $\mfa(S,z)^*J\mfa(S,z)$,
relations \eqref{I4}--\eqref{I6} yield
\begin{align}&\nn
 \mfa_{21}(z)\mfa_{22}(z)^*+\mfa_{22}(z)\mfa_{21}(z)^*
\\ &  \label{I8} 
=
\I (\ov{z}-z)\Phi_2^*(I-zA^*)^{-1}S^{-1}(I-\ov{z} A)^{-1}\Phi_2 >0 \quad \big(z\in \BC_+\big). 
\end{align} 
An important equality
\begin{align}& \label{I8+}
\mfa(S,z)J\mfa(S,\ov{z}))^*= J 
\end{align} 
follows from \eqref{I5} as well. We will need the following notation.
\begin{Nn} \label{NnPair}
A pair $\{R(z), \, Q(z)\}$  is called nonsingular, with property-$J$ if
$R(z)$ and $Q(z)$ are meromorphic $p\times p$ matrix functions in $\BC_+$ satisfying relations
 \begin{align} &  \label{I9}
R(z)^*R(z)+Q(z)^*Q(z)>0, \quad \begin{bmatrix}R(z)^* & Q(z)^* \end{bmatrix}J \begin{bmatrix}R(z) \\ Q(z) \end{bmatrix}\geq 0
 \end{align} 
$($excluding, possibly, some isolated points $z\in \BC_+)$. 
\end{Nn}

 It follows from \eqref{I8} and \eqref{I9} (see, e.g., \cite[Proposition 1.43]{SaSaR}) that 
\begin{align}& \label{I10}
\det\bigl( {\mathfrak A}_{21}(z ) R(z)+{\mathfrak A}_{22}(z)Q(z)\bigr)\not=0.
\end{align} 
Thus, the linear-fractional transformations 
\begin{equation}\label{I11}
\varphi (z) = \I \bigl( {\mathfrak A}_{11}(z)R(z)+
{\mathfrak A}_{12}(z)Q(z)\bigr) \bigl( {\mathfrak A}_{21}(z ) R(z)+{\mathfrak A}_{22}(z)Q(z)\bigr)^{-1}
\end{equation}
are well defined in the points $z\in\BC_+$ where $I-zA^*$ is invertible and \eqref{I9} holds.
\begin{Nn}\label{Nncln}
The set of matrix functions $\vp(z)$ of the form \eqref{I11}, which correspond to some fixed $S$-node
and various nonsingular pairs with proper-ty-$J$ is denoted by $\cln(\mfa(S))$.
\end{Nn}
According to \eqref{I7}, \eqref{I9} and \eqref{I11}, we have
\begin{align}& \label{I12}
\I (\vp(z)^*-\vp(z))\geq 0 \quad (z\in \BC_+).
\end{align} 
This implies that the matrix  functions $\vp(z)$ do not have singularities in $\BC_+$ (even if $\mfa(S,z)$ has poles), belong to Herglotz class and admit Herglotz representation:
\begin{align}& \label{I13}
\vp(z )=\g z +\theta+\int_{-\infty}^{\infty}\frac{1+t z }{(t-z )(1+t^2)}d\mu(t) \quad (\g\geq 0, \quad \theta=\theta^*), \\
& \label{I14}
  \int_{-\infty}^{\infty}{(1+t^2)^{-1}}{d\mu(t)}<\infty,
\end{align} 
where $\mu(t)$ is a nondecreasing $p\times p$ matrix function (and $d\mu$ is often called matrix measure).
The solutions of interpolation problems for the structured operators $S$ are   presented in \cite{SaL2} as the
triples $\{\g,\theta,\mu\}$. They are described  using linear-fractional transformations  \eqref{I11}.
The expressions 
\begin{align}& \label{I15}
\rho(z, \la):=\I (\la-z)\Phi_2^*(I-zA^*)^{-1}S^{-1}(I-\la A)^{-1}\Phi_2
\end{align} 
and $\rho(z, \ov{z})$ are important characteristics of the structured operators $S$. Therefore,
the interrelations between $\rho(z, \ov{z})$ and $\mu$, which we study in this paper, are of essential interest in interpolation
and asymptotic analysis.

In Section \ref{ApA}, the factorisation $\mu^{\prime}(t)=G_{\mu}(t)^*G_{\mu}(t)$ is considered.
Some results and references from the important survey \cite{KiKa} by V.E.~Katsnelson and B. Kirstein  were
used in that section. We note that our Theorem \ref{TmArKr} on the interrelations between $\rho(z, \ov{z})$ and $\mu$ is expressed
in terms of $G_{\mu}(z)$ (see some examples in \cite{ALS87, ALS-ConstrAppr}). The interesting
note \cite{ArKr1} by D.Z. Arov and M.G.~Krein  is applied in Section \ref{ArKr} in order to obtain this
theorem. In the spirit of \cite{ArKr1}, the inequality \eqref{B13!} in Theorem \ref{TmArKr} may be called the {\it entropy inequality}.
Finally, in Section~\ref{Out} we formulate and prove a generalisation (without various restrictions from the earlier interpolation theorem  \cite{ISa}) 
of \cite[Lemma]{ALS87}, which gives sufficient conditions for the relations \eqref{B19} to hold. Relations \eqref{B19} are the main requirements in 
Theorem \ref{TmArKr}.  We note that the proof of  \cite[Lemma]{ALS87} (generalised here)
was not published in \cite{ALS87} or later.

{\bf Notations.} The majority of the notations have been explained above. In addition, the notation $\BD$ stands for the unit disk $|\ze|<1$, the notation
${\rm tr}(Z)$ stands for the trace of $Z$, and $\Re(Z)$ denotes the real part (of either a scalar or a square  matrix $Z$).

\section{Factorisation of $\mu^{\prime}(t)$}\label{ApA}
\setcounter{equation}{0}
Factorisation of positive semi-definite integrable matrix functions
is one of the classical domains connected with the names of  A. Beurling, N. Wiener, P.R.~Masani, H. Helson, D. Lowdenslager, M.G. Krein, 
D.~Sarason and many others (see, e.g., \cite{HeLo} and the useful bibliography in \cite{KiKa}). The main results deal with the functions and matrix functions on the
unit disk. The results on the real axis follow from the classical conformal mappings between $\BC_+$ and the unit disk $\BD$.  Recall that the functions
\begin{align}& \label{A1}
z=(\overline{z_0}\zeta -z_0 )(\zeta -1)^{-1} \qquad (z_0\in \BC_+)
\end{align} 
map the unit disk onto $\BC_+$.  In spite of the simplicity of these mappings, it would be convenient to reformulate the factorisation results in terms
of the real axis (so far, we know the book \cite{Koo}, where the scalar case is treated, and some constructive results in \cite{HS} for the
matrix function case). Meanwhile, following \cite{ALS-ConstrAppr} we introduce the class $\wh H$ of matrix functions on $\BC_+$.
First recall that an outer (or maximal in the terminology of \cite{Priv}) function $h(\ze)$   is an analytic in $\BD$ function,
which admits representation
\begin{align}& \label{A2!} 
h(\ze)=\a \exp\left\{\frac{1}{2\pi}\int_{0}^{2\pi}\ln\big(\om(\vt)\big)\frac{\E^{\I \vt}+ \ze}{\E^{\I \vt}- \ze} d\vt \right\},
\end{align}
where $|\a|=1$ and $\ln\big(\om(\vt)\big)$ is integrable on $[0,2\pi]$. 
Clearly, the product of outer matrix functions is outer as well. Representation \eqref{A2!} yields that $\om(\vt)=\big|h\big(\E^{\I \vt}\big)\big|$.
(We note that functions on $\BC_+$ and $\BD$ considered in this paper have non-tangential boundary values almost everywhere
on $\BR$ and unit circle $|\ze|=1$, respectively.)
\begin{Dn}\label{wtH}  The $p \times p$ matrix function $G(z)$ belongs to $\wh H$ if 
the entries of 
\begin{align}& \label{A1+}
\wh G(\ze)=G\big((\overline{z_0}\zeta -z_0 )(\zeta -1)^{-1}\big)
\end{align}
belong to the Hardy class $H^2(\BD)$ and $ \det \left(\wh G(\zeta)\right)$ is an outer function.
\end{Dn}
\begin{Rk} \label{RkWH} We note that the accent ``widehat"  denotes $($as in \eqref{A1+}$)$ the transfer from the matrix function depending on $z\in \BC_+$ to
the matrix function depending on $\ze\in \BD$  via the substitution $z=(\overline{z_0}\zeta -z_0 )(\zeta -1)^{-1}$. 

The accent ``breve" denotes the inverse
mapping:
\begin{align}& \label{A2-}
\breve f(z)=f\big((z-z_0)/(z-\ov{z_0})\big).
\end{align} 

The corresponding limit functions on $\BR$ and on the unit circle are also denoted by  $\breve f$ and $\wh G$, respectively.
\end{Rk}
\begin{Rk} \label{RkOut} In order to  invert matrix functions we are interested in, we recall some definitions and properties of  Smirnov class $D$ of functions
and of the outer  matrix functions, which are used in this paper and  in \cite{ArKr1, KiKa} as well. 
A holomorphic in $\BD$ function $f$ belongs to $D$ if it may be represented as a ratio of a function from $H^{\infty}(\BD)$ and of an outer function from $H^{\infty}(\BD)$
$($see, e.g., \cite[p. 41]{KrSp}$)$. In particular, the function $h(z)$ given in \eqref{A2!}  and $1/h(z)$ belong $D$. If $f,g\in D$, than $fg\in D$ and $(f+g)\in D$.
 According to another $($equivalent$)$ definition of $D$ from \cite[p.~116]{Priv} and Polubarinova-Kochina theorem
\cite[p.~114]{Priv}, all the functions from the Hardy classes $H^{\delta}(\BD)$ belong to $D$. 
The notation $D^{(p\times \ell)}$ stands for the  class of $p\times \ell$ matrix functions with the entries belonging to  $D$.
A matrix function $f$ is called an outer matrix function 
if  $f\in D^{p\times p}$ and its determinant  is an outer function  \cite[p.~58]{KrSp}. From the facts above, it is clear that the inverse to the outer matrix function
exists and is again an outer matrix function.
\end{Rk} 
Next, assume that a nondecreasing  matrix function $\mu(t)$ satisfies \eqref{I14} and Szeg{\H o} condition
\begin{align}& \label{A2}
\int_{-\infty}^{\infty}(1+t^2)^{-1}\ln\big(\det\mu^{\prime}(t)\big)dt>-\infty ,
\end{align} 
where $\mu^{\prime}$ is the derivative of the absolutely
continuous part of $\mu$.
Setting
\begin{align}& \nn
t=(\overline{z_0}\zeta -z_0 )(\zeta -1)^{-1}=\overline{z_0}+(\overline{z_0}-{z_0})(\zeta-1)^{-1}, \quad \zeta=\E^{\I\theta} \quad( 0\leq \theta <2\pi),
\end{align} 
and taking into account that $\ov{\zeta}(\zeta-1)=-\ov{(\zeta-1)}$ and $\ov{\zeta}(\overline{z_0}\zeta -z_0 )=-\ov{\big(\overline{z_0}\zeta -z_0 \big)}$ we obtain
\begin{align}& \label{A3}
\frac{dt}{1+t^2}=-\frac{\I(\overline{z_0}-{z_0})\zeta d\t}{(\zeta-1)^2+(\overline{z_0}\zeta -z_0 )^2}=
\frac{\I(\overline{z_0}-{z_0})d\t}{|\zeta-1|^2+|\overline{z_0}\zeta -z_0 |^2}.
\end{align} 
From \eqref{I14} and \eqref{A3}, it follows that 
$$\int_0^{2\pi}\mu^{\prime}\left(\frac{\overline{z_0}\E^{\I \t} -z_0 }{\E^{\I \t} -1}\right)d\t<\infty.$$
In other words, $\mu^{\prime}\left(\frac{\overline{z_0}\E^{\I \t} -z_0 }{\E^{\I \t} -1}\right)$ is integrable. Since $\mu^{\prime}(t)\geq 0$ and
$\mu^{\prime}\left(\frac{\overline{z_0}\E^{\I \t} -z_0 }{\E^{\I \t} -1}\right)$ is integrable, the first two conditions of  \cite[Theorem 9]{HeLo} are fulfilled.

From  \eqref{A2} and \eqref{A3}  we derive that $\int_0^{2\pi}\ln\det\left(\mu^{\prime}\left(\frac{\overline{z_0}\E^{\I \t} -z_0 }{\E^{\I \t} -1}\right)\right)d\t>-\infty$,
and, in view of the equality $\ln(\det  (\mu^{\prime}))= {\rm tr}(\ln (\mu^{\prime}))$, the condition (74) of   \cite[Theorem 9]{HeLo} is also fulfilled.
Thus (see also \cite{KiKa}), we have the following proposition.
\begin{Pn}\label{PnFact} Let a nondecreasing  matrix function $\mu(t)$ satisfy \eqref{I14} and Szeg{\H o} condition \eqref{A2}.
Then, there is a factorisation
\begin{align}& \label{A4}
\mu^{\prime}(t)=G_{\mu}(t)^*G_{\mu}(t),
\end{align} 
where $G_{\mu}(z )\in \wh H$ and $G_{\mu}(t)$ is the boundary value function of  $G_{\mu}(z )$.
Moreover, $G_\mu(z)$  may be chosen uniquely up to a constant unitary factor from the left.
\end{Pn}
In the proposition above, we reversed the order of factors in
\cite[Theorem~9]{HeLo}, which does not matter (see, e.g., \cite[p. 195]{HeLo}).  
The  uniqueness of  $G_{\mu}(z)$ (up to a constant unitary factor from the left).
follows, for instance, from \cite[Lemma 4.1]{KiKa}. Indeed, let $G_{\mu}\in \wh H$
and $F_{\mu}\in \wh H$ satisfy \eqref{A4}. Than, $\wh G_{\mu}(\ze) \wh F_{\mu}(\ze)^{-1}$
is unitary for $|\ze|=1$. According to the ``maximum principle of V.I. Smirnov for matrix functions"
(see \cite{KiKa}), $\wh G_{\mu}(\ze) \wh F_{\mu}(\ze)^{-1}$ belongs also to the Schur class of
contractive matrix functions. Therefore, $\wh G_{\mu}(\ze) \wh F_{\mu}(\ze)^{-1}$ is an inner
matrix function. Clearly, $\wh G_{\mu}(\ze) \wh F_{\mu}(\ze)^{-1}$ is invertible and its inverse
belongs $D^{p\times p}$. Now, in view of \cite[Lemma 4.1]{KiKa},  we see that $\wh G_{\mu}(\ze) \wh F_{\mu}(\ze)^{-1}$
is a constant (and so unitary) matrix.

\begin{Rk}\label{RkMob}
The definition of $\wh H$ does not depend on the choice of $\break z_0\in \BC_+$ because $f(\ze)\in D$ yields
$$g(\eta):=f\left(\frac{a\eta-b}{\ov{b}\eta-\ov{a}} \right)\in D \quad {\mathrm{for}} \quad b\not=0, \,\, |a|>|b|, \,\, |\eta|<1.$$
\end{Rk}
\section{Arov-Krein results: reformulation}\label{ArKr}
\setcounter{equation}{0}
Interesting inequalities (and asymptotic results) based on linear fractional transformations have been
published by D.Z. Arov and M.G. Krein in \cite{ArKr1} (see also some related proofs in \cite{ArKr2}).
In order to use them in our work, certain reformulations are necessary.

The $2p \times 2p$ meromorphic matrix functions $\cla(\ze)$ ($|\ze |< 1$), such that
\begin{align}& \label{B1}
\cla(\ze)^*J\cla(\ze)\geq -j , \quad j:=\begin{bmatrix}   I_p & 0 \\ 0 & -I_p \end{bmatrix},
\end{align} 
have been considered in \cite{ArKr1}. Note that changes in this section  of  several notations from \cite{ArKr1}  were caused by the
notations used in the previous sections. For instance,    the notation $A(z)$ have been used in \cite{ArKr1} instead
of $\cla(\ze)$ and $n$ have been used instead of $p$.

Let us assume that an $S$-node $\{A,S,\Pi\}$ is given, that is, the identity \eqref{I1} holds, $A,\, S \in \clb(\clh)$, and $S=S^*$.  Assume that $S$ satisfies the first relation in \eqref{I5} and that 
the operators $I-zA^*$ have  bounded inverses for $z\in \BC_+$ (excluding, possibly, some isolated points, which are the poles of the
matrix function $\mfa(S,z)$ given by  \eqref{I2}).

It is easy to see that $J$ and $j$ are unitarily similar:
\begin{align}& \label{c13}
J=KjK^*, \quad {\mathrm{where}} \quad K:= \frac{1}{\sqrt{2}}\begin{bmatrix}I_p & -I_p \\ I_p & I_p \end{bmatrix}, \quad K^*=K^{-1}.
\end{align}
Therefore, we also have $J=KJ(-j)JK^*$ and so \eqref{I7} implies  that the matrix function
\begin{align}& \label{B3}
\cla(\ze):=\wh \mfa(S,\ze)KJ,
\end{align} 
where $\wh \mfa$ is defined in Remark \ref{RkWH}, satisfies \eqref{B1}. The pairs $\{q(z),I_p\}$ $\break (z\in \BC_+)$ such that
$q$ are $p\times p$ holomorphic matrix functions and
\begin{align}& \label{B3+}
q(z)^*q(z)\leq I_p,
\end{align} 
are nonsingular, with property-$(-j)$ and all nonsingular pairs with property-$(-j)$ have the form $ \{q(z)a( z), a(z)\}$, where $\det \big(a(z)\big)\not=0$.
Since $\break JK^*JKJ=-j$, the transformation
\begin{align}& \label{B5}
\begin{bmatrix}R(z) \\ Q(z)\end{bmatrix}=KJ\begin{bmatrix}q(z)a( z) \\ a( z)\end{bmatrix}
\end{align} 
maps the set of nonsingular pairs with property-$(-j)$ onto the set of nonsingular pairs with property-$J$. 

Next, we compare linear fractional transformations \eqref{I11} with  linear fractional transformations from \cite{ArKr1}:
\begin{align}& \label{B4}
f(\ze)=\bigl( {\cla}_{11}(\ze)\wh q(\ze)+
{\cla}_{12}(\ze)\bigr) \bigl( {\cla}_{21}(\ze ) \wh q(\ze)+{\cla}_{22}(\ze)\bigr)^{-1},
\end{align} 
where $\cla_{ij}$ are $p\times p$ blocks of $\cla$ and $\wh q(\ze)$ are contractive in $\BD$. 
In view of \eqref{B3}, \eqref{B5},
the set $\cln(\mfa(S))$ of matrix functions $\vp(z)$ given by \eqref{I11} and the set  $\cln(\cla)$ of  matrix functions $f(\ze)$ given by \eqref{B4}
are connected by the relation
\begin{align}& \label{B6}
\wh \vp(\ze)=\I f(\ze) \quad {\mathrm{or}} \quad \vp(z)=\I \breve f(z),
\end{align} 
where $\breve f$ is defined in \eqref{A2-}.
It follows from the well-known Stieltjes-Perron formula (see, e.g., \cite[Theorem~2.2~(v)]{GT} and the bibliography there) that for the boundary values
of $\vp$ and $f$ considered in \eqref{B6} we have
\begin{align}& \label{B7}
\Im (\vp(t))=\Re\big(\breve f(t)\big)=\pi \mu^{\prime}(t) \quad (t\in \BR),
\end{align} 
where $\mu$ is the nondecreasing matrix function in the Herglotz representation \eqref{I13} of  $\vp$.

Now, let us consider $\chi(\ze):=-\cla_{22}(\ze)^{-1}\cla_{21}(\ze)$. In view of \eqref{B3}, we have
\begin{align}& \label{B8}
\chi(\ze)=(\wh\mfa_{21}( \ze)+\wh \mfa_{22}(\ze))^{-1}(\wh \mfa_{21}( \ze)-\wh \mfa_{22}( \ze)),
\end{align}
where $\wh \mfa_{ik}(\ze)=\wh \mfa_{ik}(S,\ze)$ are $p\times p$ blocks of $\wh \mfa(S,\ze)$. 
Clearly, the requirement 
\begin{align}& \label{B11+}
\|\chi(\ze)\|<1 \quad {\mathrm{for}} \quad z\in \BD
\end{align} 
 in \cite[Theorems 1 and 3]{ArKr1} is equivalent to the condition $\chi(\ze)\chi(\ze)^*<I_p$. Therefore, using \eqref{B8}, after easy transformations  
we derive that \eqref{B11+} is equivalent to 
\begin{align}& \label{B9}
\wh \mfa_{21}( \ze)\wh \mfa_{22}(\ze)^*+\wh \mfa_{22}( \ze)\wh \mfa_{21}(\ze)^*>0.
\end{align} 
Under conditions \eqref{I5}, the inequality \eqref{B9} is immediate from \eqref{I8}.  Thus, the requirement \eqref{B11+}in \cite{ArKr1}
is fulfilled.

By virtue  of \eqref{B1} and  \eqref{B11+}, the conditions of  \cite[Theorem 1]{ArKr1} are satisfied. 
Hence, the conditions of \cite[Theorem 3]{ArKr1} are fulfilled
for the case where the function
\begin{align}& \label{B14}
\cla_{22}(\ze)=\frac{1}{\sqrt{2}}\big(\wh\mfa_{21}( \ze)+\wh \mfa_{22}(\ze)\big)
\end{align} 
is an outer matrix function. Moreover, taking into account \eqref{B3}, in that case we have
\begin{align}\nn
\Delta(\ze): &=\big(\cla_{22}(\ze)\cla_{22}(\ze)^*-\cla_{21}(\ze)\cla_{21}(\ze)^*\big)^{-1}
\\ & \label{B15}
=\big(\wh \mfa_{21}( \ze)\wh \mfa_{22}(\ze)^*+\wh \mfa_{22}( \ze)\wh \mfa_{21}(\ze)^*\big)^{-1},
\end{align} 
where $\Delta$ is an important notion from \cite{ArKr1}. It follows from \eqref{I8}, \eqref{I15} and \eqref{B15} that
\begin{align}& \label{B15+}
\Delta(\ze)=\wh \rho\big(\ze, 1/ \, \ov{\ze}\big)^{-1} \quad \big(\wh \rho\big(\ze, 1/ \, \ov{\ze}\big)>0\big).
\end{align} 

The requirement above that $\cla_{22}(\ze)$ is an outer matrix function may be rewritten in the form
\begin{align}& \label{B16}
\cla_{22}(\ze), \, \cla_{22}(\ze)^{-1} \in D^{(p\times p)},
\end{align}
where $D^{(p\times \ell)}$ is 
introduced in Remark \ref{RkOut}. Using \eqref{B14}, we rewrite $\cla_{22}(\ze)$ as
\begin{align}& \label{B17}
\cla_{22}(\ze)=\frac{1}{\sqrt{2}}\wh\mfa_{21}( \ze)\big(I_p+\wh\mfa_{21}( \ze)^{-1}\wh \mfa_{22}(\ze)\big).
\end{align}
Note that \eqref{B9} yields the invertibility of $\wh\mfa_{21}( \ze)$ and, moreover, we have
$\Re\big(\wh\mfa_{21}( \ze)^{-1}\wh \mfa_{22}(\ze)\big)>0$. Now, it is immediate that
\begin{align}& \label{B18}
\Re\big(I_p+\wh\mfa_{21}( \ze)^{-1}\wh \mfa_{22}(\ze)\big)>0, \quad \Re\Big(\big(I_p+\wh\mfa_{21}( \ze)^{-1}\wh \mfa_{22}(\ze)\big)^{-1}\Big)>0.
\end{align}
Thus, according to Smirnov's theorem (see, e.g., \cite[p. 93]{Priv}) the functions 
\begin{align}& \label{B18+}
h^*\big(I_p+\wh\mfa_{21}( \ze)^{-1}\wh \mfa_{22}(\ze)\big)h, \quad h^*\big(I_p+\wh\mfa_{21}( \ze)^{-1}\wh \mfa_{22}(\ze)\big)^{-1}h
\end{align}
belong to the classes $H^{\delta}$ for any $h\in \BC^p$ and $0<\delta<1$. Therefore, we have
$$\big(I_p+\wh\mfa_{21}( \ze)^{-1}\wh \mfa_{22}(\ze)\big), \, \big(I_p+\wh\mfa_{21}( \ze)^{-1}\wh \mfa_{22}(\ze)\big)^{-1} \in D^{p\times p},$$
and so it suffices that
\begin{align}& \label{B19}
\wh \mfa_{21}(\ze), \, \wh \mfa_{21}(\ze)^{-1} \in D^{(p\times p)}
\end{align}
for \eqref{B16} to hold. 

Summing up the reformulated conditions of  \cite[Theorem 3]{ArKr1}, we state the result.
\begin{Tm}\label{TmArKr} Let an $S$-node $\{A,S, \Pi\}$, such that 
the operators $I-zA^*$ have  bounded inverses for $z\in \BC_+$ $($excluding, possibly, some isolated points, which are the poles of the
matrix function $\mfa(S,z))$, be given and let relations \eqref{I5} and \eqref{B19} hold. Assume that $\vp(z)\in \cln(\mfa(S))$ and the
matrix function $\mu(t)$ in Herglotz representation \eqref{I13} of $\vp$ satisfies Szeg{\H o} condition \eqref{A2}. 

Then, we have
\begin{align}& \label{B13!}
2 \pi  G_{\mu}(z)^* G_{\mu}(z)\leq \rho\big(z, \ov{z}\big)^{-1}  \quad (z \in \BC_+).
\end{align} 
The equality in \eqref{B13!} holds in some point $z=\la \in \BC_+$ if and only if our $\vp(z)\in \cln(\mfa(S))$ is generated $($see \eqref{I11}$)$ by the
constant pair $\{R,Q\}:$
\begin{align}& \label{B31}
R(z)\equiv \mfa_{22}(S,\la)^*, \quad Q(z)\equiv \mfa_{21}(S,\la)^*.
\end{align} 
\end{Tm}
\begin{proof}
First note that $\wh \mfa(S,\ze) KJ$ satisfies conditions of   \cite[Theorem 3]{ArKr1} (as shown above).
According to \eqref{A4} and \eqref{B7}, we have
\begin{align}& \label{B12}
\pi \wh G_{\mu}(\ze)^*\wh G_{\mu}(\ze)=\Re\big( f(\ze)\big) \quad {\mathrm{for}} \quad |\ze |=1.
\end{align} 
Comparing \eqref{B12} with the factorisation at the beginning of paragraph 3 in \cite{ArKr1}, we see
that $\sqrt{\pi}\wh G_{\mu}(\ze)$ coincides with $\vp_f(\ze)$ in the notations of \cite{ArKr1}. Thus,
the inequality in \cite[Theorem 3]{ArKr1} takes the form
\begin{align}& \label{B13}
2 \pi \wh G_{\mu}(\ze)^*\wh G_{\mu}(\ze)\leq \Delta(\ze) \quad (\ze \in \BD),
\end{align} 
where the equality is attained at some $\ze =\wt \ze$  if and only if $\wh q(\ze)\equiv \chi(\wt \ze)^*$  (for $\wh q$ in \eqref{B4}).
  (We note that \cite{ArKr1} contains a misprint and the coefficient
``2" from \eqref{B13} is missing there.) 
Recall that the meaning of the accent ``widehat" in the paper is explained in Remark \ref{RkWH}.
Hence, taking into account \eqref{B15+}, we see that the inequality \eqref{B13} is equivalent to \eqref{B13!}.
The interrelations between $\wh q$ and  pairs $\{R,Q\}$ in \eqref{I11}, which generate $\vp(z)$ corresponding
$f(\ze)$ (as in \eqref{B7}), is given by \eqref{B5}. Therefore, condition $\wh q(\ze)\equiv \chi(\wt \ze)^*$, where
$\chi(\ze)$ has the form \eqref{B8}, transforms into \eqref{B31}.

\end{proof}

\section{Outer matrix function}\label{Out}
\setcounter{equation}{0}
Relation \eqref{B19} on $\wh \mfa_{21}(\ze)$ is the most complicated to check among the requirements stated in Theorem \ref{TmArKr}.
Here, we present some sufficient conditions for \eqref{B19} to hold.
It follows from \eqref{I2} that 
\begin{align}& \label{C0}
\mfa(S,z)=I_{2p}-\I z\Pi^*(I-zA^*)^{-1}S^{-1}\Pi J.
\end{align}
In particular, $ \mfa_{21}(S,z)$ has the form
\begin{align}& \label{C1}
c(z):=\mfa_{21}(S,z)=-\I z\Phi_2^*(I-zA^*)^{-1}S^{-1}\Phi_2.
\end{align} 
Proposition \ref{PnPM} on $\wh c(\ze)=\wh \mfa_{21}(\ze)$ is helpful in the applications of  Theorem~\ref{TmArKr}. 
\begin{Pn}\label{PnPM}  Let a triple $\{A,S,\Pi\}$ form an $S$-node and let
the operators $I-zA$ have bounded inverses for $z$ in the domains  $\{z: \, \Im(z)\leq 0\}$ and $\{z: \, \Im(z)>0, \,\, |z| \geq r_0\}$ for some $r_0>0$.
Assume that  relations \eqref{I5}  hold. Finally, let
\begin{align}& \label{Z2}
{\ov{\lim}}_{r\to \infty}\big(\ln(\clm(r))\big/r^{\vk}\big)<\infty  
\end{align} 
for some  $0<\vk<1$ and $\clm(r)$ given by
\begin{align}& \label{Z3}
\clm(r)=\sup_{r_0<|z|<r} \|(I-z A)^{-1}\|.
\end{align} 

Then, we have $\wh c(\ze), \, \wh c(\ze)^{-1} \in D^{(p\times p)}$, that is, $\wh c(\ze)$ is an outer matrix function.
\end{Pn}

The proof of Proposition \ref{PnPM} consists of the proofs of two lemmas below.
\begin{La}\label{LaPlus1} Let a triple $\{A,S,\Pi\}$ form an $S$-node, let the operator $S$ have a bounded inverse, 
and let the operators $I-z A^*$ have bounded inverses for $\Im(z)\geq 0$.
Assume that 
\begin{align}& \label{C5}
{\ov{\lim}}_{r\to \infty}\big(\ln(M(r))\big/r^{\vk}\big)<\infty \quad {\mathrm{for\,\, some}} \quad  0<\vk<1,
\end{align} 
where 
\begin{align}& \label{C6}
M(r)=\sup_{|z|<r, \,\, \Im(z)\geq 0} \|(I-z A^*)^{-1}\|.
\end{align} 
Then, $\wh c(\ze)\in  D^{(p\times p)}$.
\end{La}
\begin{proof} 
It follows from \eqref{C1} that
\begin{align}& \label{C7}
\|c(z)\| \leq \clc_1 \, |z| \, \|(I-zA^*)^{-1}\| \quad {\mathrm{for \,\, some}} \quad \clc_1>0.
\end{align} 
According to Remark \ref{RkMob}, it suffices to prove our lemma for the case $z_0=\I$ in \eqref{A1}. Thus, we set
\begin{align}& \label{C8}
z=- \I (\ze+1)(\ze -1)
\end{align} 
Clearly, $\wh c(\ze)$ is analytic on the closed unit disk $|\ze| \leq1$  excluding, possibly,
the point $\ze =1$.  By virtue of \eqref{C5}--\eqref{C7}, we have
\begin{align}& \label{C9}
|\wh c_{ij}(\ze)|\leq \clc_2 |\ze-1|^{-1}\E^{\clc_3 |\ze -1|^{-\vk}}
\end{align} 
for some $\clc_2,\clc_3$ and the entries $\wh c_{ij}$ of $\wh c$. We also note that
\begin{align}& \label{C9+}
2\big|r\E^{\I  \vt}-1\big|>\big|\E^{\I  \vt}-1\big| \quad (0<|\vt|<\pi, \quad 0<r<1).
\end{align} 
In order to prove that
$\wh c(\ze)\in D^{(p\times p)}$, one needs to show (in the disk $|\ze|\leq 1$) the equalities
\begin{align}& \label{C10}
\lim_{r \to 1}\int_{-\pi}^{\pi}\ln^{+}\big|\wh c_{ij}\big(r\E^{\I \vt}\big)\big|d\vt=\int_{-\pi}^{\pi}\ln^{+}\big|\wh c_{ij}\big(\E^{\I \vt}\big)\big|d\vt
\end{align}
for $1\leq i,j \leq p$, where $\ln^{+}a$ equals $\ln a$ for $a \geq1$ and equals $0$ for $a<1$. Equalities \eqref{C10} follow from \eqref{C9}, \eqref{C9+}
and the continuity of $\wh c(\ze)$ (excluding, possibly, the point $\ze =1$) if we split the domain $[-\pi,\pi]$ of integration in \eqref{C10} into
the domains $[-\pi,\delta]\cup [\delta, \pi]$ and $(-\delta, \delta)$ with $\delta$ tending to zero. In view of one of the equivalent
definitions of the class $D$ (see \cite[p. 116]{Priv}), equalities \eqref{C10} imply $\wh c(\ze)\in  D^{(p\times p)}$.
\end{proof}
\begin{Rk} Similar to \eqref{C9+}, one obtains a slightly more general relation, which we will need further in the text. That is, for $\vt,\vt_k\in \BR, \,\, |\vt -\vt_k|<\pi$
we have
\begin{align}& \label{C9!}
2\big|r\exp\{\I  \vt\}-\exp\{\I  \vt_k\}\big|>\big|\exp\{\I  \vt\}-\exp\{\I  \vt_k\}\big| .
\end{align} 
\end{Rk}
\begin{Rk}  We note that the requirement of the invertibility of $(I-z A^*)$ for $\Im(z)\geq 0$ does not hold for our standard triple
$\{A,S,\Pi\}$ corresponding to Toeplitz matrices $S$ $($see \cite[(2.5)--(2.8)]{ALS-ConstrAppr}$)$. However, this triple is easily substituted by the triple
$\{\wt A,S, \wt\Pi\}$, where $\wt A=-A$ and $\wt \Pi=\begin{bmatrix} -\Phi_2 & \Phi_1\end{bmatrix}$ and the
corresponding asymptotic result is obtained in this way in \cite[Example~1]{ALS87}.
\end{Rk}
Let the conditions of Lemma \ref{LaPlus1} and relations \eqref{I5} be satisfied. Then,
we easily see that $c(z)$ is analytic in the domains 
\begin{align}& \label{Z1}
\{z: \, -\wt r \leq \Re(z)\leq \wt r, \,\, \Im(z)\geq -\ve(\wt r)\}
\end{align}
 for all the values $\wt r >0$ and some values $\ve(\wt r)>0$.
Since $S^{-1}>0$ and $\ker \, \Phi_2=~0$, similar to \eqref{I8}, we obtain
\begin{align}& \label{I8-}
 \mfa_{21}(z)\mfa_{22}(z)^*+\mfa_{22}(z)\mfa_{21}(z)^*<0 \quad {\mathrm{for}}\quad
\Im(z) <0.
\end{align} 
Formulas \eqref{I8} and   \eqref{I8-} imply that $c(z)=\mfa_{21}(S,z)$ is invertible in these domains excluding possibly
some values $z\in \BR$. This yields that $\det \big(c(z)\big)$ may have only finite number of zeros of finite order in  each domain \eqref{Z1}
and these zeros are placed on the interval $ -\wt r<z<\wt r$. 

We again assume the correspondence \eqref{C8} between $z$ and $\ze$. Recall that $\wh c(\ze)\in D^{(p\times p)}$
and so $\det \wh c(\ze)\in D$ and the minors of $\wh c(\ze)$ belong $D$ as well.  Since the functions from $D$ may be represented
as ratios of bounded functions analytic in the unit disk, the entries of $\wh c(\ze)^{-1}$ may be represented as such ratios as well.
Moreover $\wh c(z)^{-1}$ is analytic in the unit disk because $c(z)$ is invertible in $\BC_+$. Thus, the entries of $\wh c(\ze)^{-1}$
have characteristic properties of the functions from class ${\bf A}$ (see \cite[p. 82]{Priv}), that is, $\wh c(\ze)^{-1}\in {\bf A}^{(p\times p)}$.
For the entries $\big(\wh c(\ze)^{-1}\big)_{ij}$ of $\wh c(\ze)^{-1}$, it follows (see \cite[p. 78, f-la (1.1:1)]{Priv}) that 
\begin{align}& \label{C3}
\lim_{r\to 1}\int_{-\pi}^{\pi}\ln^{+}\big|\big(\wh c(r\E^{\I \vt})^{-1}\big)_{ij}\big|d\vt<\infty ,
\end{align} 
and (see \cite[p. 83, f-la (2.1:2)]{Priv}) that
\begin{align}& \label{C3+}
\int_{-\pi}^{\pi}\ln^{+}\big|\big(\wh c(\E^{\I \vt})^{-1}\big)_{ij}\big|d\vt<\infty .
\end{align} 
However, relations \eqref{C3} and \eqref{C3+} are insufficient for our purposes.
In order to show that
\begin{align}& \label{C4}
\lim_{r\to 1}\int_{-\pi}^{\pi}\ln^{+}\big|\big(\wh c(r\E^{\I \vt})^{-1}\big)_{ij}\big|d\vt=\int_{-\pi}^{\pi}\ln^{+}\big|\big(\wh c(\E^{\I \vt})^{-1}\big)_{ij}\big|d\vt ,
\end{align} 
that is, $ \wh c(\ze)^{-1} \in D^{(p\times p)}$, we need stronger conditions.
\begin{La}\label{LaMinInD}
 Let a triple $\{A,S,\Pi\}$ form an $S$-node and let
the operators $I-zA^*$ have bounded inverses for $z$ in the domain  $\Im(z)\geq 0$ as well as for $z$ in the domain $\{z: \, \Im(z)<0, \,\, |z| \geq r_0\}$ for some $r_0>0$.
Assume that  relations \eqref{I5} and \eqref{C5} hold. Finally, let
\begin{align}& \label{C2'}
{\ov{\lim}}_{r\to \infty}\big(\ln(\wt M(r))\big/r^{\vk}\big)<\infty \quad {\mathrm{for\,\, some}} \quad  0<\vk<1,
\end{align} 
where 
\begin{align}& \label{C3'}
\wt M(r)=\sup_{\Im(z)\leq 0, \, r_0<|z|<r} \|(I-z A^*)^{-1}\|.
\end{align} 

Then, we have $\wh c(\ze), \, \wh c(\ze)^{-1} \in D^{(p\times p)}$.
\end{La}
\begin{proof} 
Since the conditions of Lemma \ref{LaPlus1} are satisfied, we see that $\break \wh c(\ze) \in D^{(p\times p)}$.
In order to show that $\wh c(\ze)^{-1} \in D^{(p\times p)}$, we need some preparations.

Inequality \eqref{I7} yields
\begin{align}& \label{C11}
\begin{bmatrix} a(z)^* & c(z)^*\end{bmatrix}J \begin{bmatrix} a(z) \\ c(z)\end{bmatrix}\geq 0 \quad \big(a(z):=\mfa_{11}(S,z)\big).
\end{align} 
Formula \eqref{C11} implies that $a(z)c(z)^{-1}+ \big(c(z)^{-1}\big)^*a(z)^*\geq 0$.
Now, using Smirnov's theorem similar to the considerations for the functions \eqref{B18+}, we obtain 
\begin{align}& \label{C12}
\Up(\ze)=\{\Up_{i,j}\}_{i,j=1}^p(\ze):=\wh a(\ze) \wh c(\ze)^{-1}\in D^{(p \times p)}.
\end{align} 

By virtue of \eqref{I6}, we have
\begin{align}& \label{C13}
\mfa(S,z)^{-1}=J\mfa(S,\ov{z})^*J.
\end{align} 
Recall that $a=\mfa_{11}$ and $c=\mfa_{21}$. We will use the following representation:
\begin{equation}  \label{C14}
c(z)^{-1}=\begin{bmatrix} I_p & 0\end{bmatrix}\mfa(S,z)^{-1}\begin{bmatrix} a(z) \\ c(z)\end{bmatrix}c(z)^{-1}=\clv(z)
\begin{bmatrix} a(z)c(z)^{-1} \\ I_p\end{bmatrix},
\end{equation} 
where 
\begin{align}& \label{C15}
\clv(z):=\begin{bmatrix} I_p & 0\end{bmatrix}\mfa(S,z)^{-1}=\begin{bmatrix} I_p & 0\end{bmatrix}J\mfa(S,\ov{z})^*J.
\end{align}
 In view of \eqref{C12} and \eqref{C14}, it suffices to show that $\wh \clv(\ze)$ has a good behaviour near the unit circle
in order to prove this proposition. 

Indeed, it follows from \eqref{I2} and \eqref{C15} that
\begin{align}& \label{C16}
\clv(z)=\begin{bmatrix} I_p & 0\end{bmatrix}+ \I z \Phi_1^*S^{-1}(I-zA)^{-1}\Pi J.
\end{align}
Similar to \eqref{C9}, one can show that
\begin{align}& \label{C17}
|\wh \clv_{ij}(\ze)|\leq \clc_4 |\ze-1|^{-1}\E^{\clc_5 |\ze -1|^{-\vk}} \quad (1\leq i\leq p, \,\, 1\leq j \leq 2p)
\end{align} 
in some open neighbourhood of $\ze=1$ in the disk $|\ze| \leq 1$.  Let us fix  a set 
$$\Om=\Om(\wt r, \wt \de)=\{(\de,r): \, 0<\de \leq \wt \de, \,\, 0<\wt r\leq r\leq1\},$$
where \eqref{C17} holds for all $\ze= r\E^{\pm \I  \de}$ such that $(\de,r)\in \Om $.
From \eqref{C9+}, \eqref{C14} and \eqref{C17}
we derive  (for $(\de,r)\in \Om $) that
\begin{align}\nn
\int_{-\de}^{\de}\ln^{+}\big|\big(\wh c(r\E^{\I \vt})^{-1}\big)_{i,j}\big|d\vt \leq & \int_{-\de}^{\de}\ln^{+}\big( \clc_6 \big|\E^{\I \vt}-1\big|^{-1}\exp\{\clc_7 \big|\E^{\I \vt} -1|^{-\vk}\}\big)d\vt
\\ & \label{C18}
+\sum_{k=1}^p\int_{-\de}^{\de}\ln^{+}\big|\Up_{kj}\big(r\E^{\I \vt}\big)\big|d\vt+ 2\de\ln(p+1).
\end{align}

Let us consider  $\Up(\ze)$ given in \eqref{C12}. It follows from \eqref{C0} that $\wh a(\ze)$ and $\wh c(\ze)$ are holomorphic on the disk $|\ze | \leq 1$ excluding, possibly, the point $\ze=1$.
Using our  considerations on   $\det\big(c(z)\big)$ before this lemma, we see that  $\wh c(\ze)^{-1}$
 is holomorphic on the disk $|\ze | \leq 1$ excluding, possibly, the point $\ze=1$ and poles on the unit circle. Therefore, $\Up(\ze)$
 is also holomorphic on the disk $|\ze | \leq 1$ excluding, possibly, the point $\ze=1$ and poles on the unit circle.
Moreover, there is a finite number
of these poles $\ze=\exp\{\I \vt_{\ell}\}$ on the arc $\ze=\E^{\I \vt}$, where $|\vt | \geq \de$. Hence, using \eqref{C9!}, we see that
\begin{align}& \label{C19}
\int_{-\pi}^{\de}\ln^{+}\big|\Up_{kj}\big(r\E^{\I \vt}\big)\big|d\vt \to \int_{-\pi}^{\de}\ln^{+}\big|\Up_{kj}\big(\E^{\I \vt}\big)\big|d\vt , \\
& \label{C20}
\int_{\de}^{\pi}\ln^{+}\big|\Up_{kj}\big(r\E^{\I \vt}\big)\big|d\vt \to \int_{\de}^{\pi}\ln^{+}\big|\Up_{kj}\big(\E^{\I \vt}\big)\big|d\vt ,
\end{align}
for $1\leq k \leq p, \,\, r\to 1$. We note that \eqref{C19} and \eqref{C20} are derived by considering separately the intervals of integration
$(\vt_{\ell}-\de_{\ell},\vt_{\ell}+\de_{\ell})$ ($\de_{\ell}\to 0$) in order to exclude the poles.
Since $\Up_{kj}(\ze)\in D$, formulas \eqref{C19} and \eqref{C20} imply a similar relation on the arc $|\vt|\leq \de$:
\begin{align}& \label{C21}
\int_{-\de}^{\de}\ln^{+}\big|\Up_{kj}\big(r\E^{\I \vt}\big)\big|d\vt \to \int_{-\de}^{\de}\ln^{+}\big|\Up_{kj}\big(\E^{\I \vt}\big)\big|d\vt .
\end{align}
It follows from \eqref{C18} and \eqref{C21} that for any $\ve >0$ there are a sufficiently small $\de$ ($\de \leq \wt \de$) and a value $r_{\ve}$ $(\wt r \leq r_{\ve}< 1)$
such that for
$r$ satisfying $r_{\ve}\leq r\leq 1$ we have
\begin{align}& \label{C22}
\int_{-\de}^{\de}\ln^{+}\big|\big(\wh c(r\E^{\I \vt})^{-1}\big)_{i,j}\big|d\vt< \ve/4 .
\end{align}
Similar to the proof of \eqref{C19} and \eqref{C20}, we obtain that (for the fixed $\de$ discussed above, some $\wt r_{\ve}$ $(\wt r \leq \wt r_{\ve}<1)$, and all $r$
such that $\wt r_{\ve} \leq r <1$)
the inequality
\begin{equation} \label{C23}
\Big| \int_{\de\leq |\vt|\leq \pi}\ln^{+}\big|\big(\wh c(\E^{\I \vt})^{-1}\big)_{i,j}\big|d\vt-\int_{\de\leq |\vt|\leq \pi}\ln^{+}\big|\big(\wh c(r\E^{\I \vt})^{-1}\big)_{i,j}\big|d\vt\Big|<\ve/2
\end{equation}
holds. Taking into account that some $\de$,  $ r_{\ve}$ and $\wt r_{\ve}$ exist for each $\ve$, we see that the inequalities \eqref{C22} and \eqref{C23}
imply \eqref{C4}, that is, $\big(\wh c(\ze)^{-1}\big)_{i,j}\in D$.  Thus, we have $\wh c(\ze)^{-1}\in D^{(p\times p)}$.
\end{proof}
Proposition \ref{PnPM} follows from Lemma \ref{LaMinInD} after minor changes including formulation of the conditions in terms of $I-zA$
instead of $I-zA^*$.



\begin{thebibliography}{AGKS}
\bibitem{ArKr1} 
Arov, D.Z., Krein, M.G.: Problem of search of the minimum of entropy in indeterminate extension problems (Russian).   {Funct. Anal. Appl.} {\bf 15}, 123--126 (1981)
\bibitem{ArKr2} 
Arov, D.Z., Krein, M.G.: Calculation of entropy functionals and their minima in indeterminate continuation problems (Russian).  Acta Sci. Math. {\bf 45}(1-4), 33--50 (1983)
\bibitem{GT}
Gesztesy, F., Tsekanovskii, E.: On Matrix-Valued Herglotz Functions. Math. Nachr. {\bf 218}, 61--138 (2000)
\bibitem{HeLo}
Helson, H., Lowdenslager, D.: 
Prediction theory and Fourier series in several variables.
{Acta Math.} {\bf 99},  165--202 (1958)
\bibitem{ISa}
Ivanchenko T.S., Sakhnovich, L. A.: An operator approach to the study of interpolation problems.
Manuscript No. 701, UK-85, deposited at the Ukrainian NIINTI (1985)
 \bibitem{KiKa}
Katsnelson, V.E., Kirstein, B.: 
On the theory of matrix-valued functions belonging to the Smirnov class.  
In: Oper. Theory Adv. Appl.  {\bf 95},  299--350, Birkh\"auser, Basel (1997) 
\bibitem{Koo}
Koosis, P.:  Introduction to $H_p$ spaces. With an appendix on Wolff's proof of the corona theorem.  Cambridge University Press, Cambridge-New York, 1980.
\bibitem{KrSp}
Krein, M.G.,  Spitkovsky, I.M.: Factorization of $\a$-sectorial matrix-valued functions on the unit circle (Russian). Mat. Issled., No. 47, 41--63 (1978)
\bibitem{Priv} 
Privalov, I.I.:  Boundary properties of analytic functions (Russian).
2d ed., Gosudarstv. Izdat. Tehn.-Teor. Lit., Moscow (1950)
\bibitem{ALS87} Sakhnovich, A.L.:
On a class of extremal problems. Math. USSR-Izv. {\bf  30}, 411--418 (1988)
\bibitem{ALS-ConstrAppr}
Sakhnovich, A.L.:
Discrete self-adjoint Dirac systems: asymptotic relations, Weyl functions and Toeplitz matrices. Constr. Approx. {\bf 55}, 641--659 (2022)
\bibitem{SaSaR}
Sakhnovich, A.L.,  Sakhnovich, L.A.,  Roitberg, I.Ya.:   Inverse Problems and Nonlinear Evolution Equations. 
 Solutions, Darboux Matrices and Weyl--Titchmarsh Functions. De Gruyter, Berlin (2013)
\bibitem{SaL1}
Sakhnovich, L.A.:
On  the  factorization  of  the  transfer
matrix function. Sov. Math. Dokl.    {\bf 17},  203--207 (1976)
\bibitem{SaL2}
Sakhnovich, L.A.:
Interpolation theory and its applications. Kluwer, Dordrecht (1997)
\bibitem{HS}
Salehi, H.:  A factorization algorithm for $q \times q$ matrix-valued functions on the real line $\BR$. Trans. Amer. Math. Soc. {\bf 124}, 468--479 (1966)
\end{thebibliography}
\end{document}